\documentclass[review]{elsarticle}
\usepackage{hyperref}
\usepackage{amsfonts}

\begin{document}

\newtheorem{thm}{Theorem}
\newtheorem{lem}[thm]{Lemma}
\newtheorem{cor}[thm]{Corollary}
\newtheorem{pro}[thm]{Proposition}
\newdefinition{dfn}[thm]{Definition}
\newproof{prf}{Proof}

%\begin{frontmatter}

\title{Two Forbidden Induced Minor Theorems for Antimatroids}
\author[cja]{Christian Joseph Altomare}
\ead{altomare@math.ohio-state.edu}
\address{The Ohio State University, 231, West 18th Avenue,
Columbus, Ohio, United States}

\begin{abstract}
Antimatroids were discovered by Dilworth in the context of lattices [4] and introduced
by Edelman and Jamison as convex geometries in[5].
The author of the current paper independently discovered (possibly infinite) antimatroids in the context of proof systems
in mathematical logic [1].
Carlson, a logician, makes implicit use of this view of proof systems as possibly infinite antimatroids in [2].
Though antimatroids are in a sense dual to matroids, far fewer antimatroid forbidden minor theorems are
known.
Some results of this form are proved in [6], [7], [8], and [9].
This paper proves two forbidden induced minor theorems for
these objects, which we think of as proof systems.

Our first main theorem gives a new proof of the forbidden induced
minor characterization of partial orders as proof systems, proved in [8]
in the finite case and
stated in [10] for what we call strong aut descendable proof systems.
It essentially states that, pathologies aside,
there is a certain unique simplest nonposet. Our second main theorem
states the new result that, pathologies aside, there is a certain unique simplest proof system
containing points $x$ and $y$ such that $x$ needs $y$ in one context,
yet $y$ needs $x$ in another.

\end{abstract}

\maketitle

%\end{frontmatter}


\section{Introduction}

\def\lad{\{}
\def\rad{\}}
\def\st{such that }
\def\cmp{comparable }
\def\icmp{incomparable }
\def\aut{autonomous }
\def\str{lengthen }
\def\strs{lengthens }
\def\strd{lengthened }
\def\strg{lengthening }
\def\wea{fatten }
\def\weas{fattens }
\def\weag{fattening }
\def\wead{fattened }
\def\wfd{well founded }
\def\wo{well order }
\def\wod{well ordered }
\def\wos{well orders }
\def\dc{downward closed }
\def\uc{upward closed }
\def\z{partial order }
\def\zs{partial orders }
\def\zd{partially ordered }
\def\pr{predecessor }
\def\prs{predecessors }
\def\pa#1#2{$(#1,#2)$}
\def\epa#1#2{(#1,#2)}
\def\part#1#2{$(#1_{#2},\le_{#2})$}
\def\epart#1#2{(#1_{#2},\le_{#2})}
\def\re#1{$\le_{#1}$}
\def\res#1{$<_{#1}$}
\def\ere#1{\le_{#1}}
\def\eres#1{<_{#1}}
\def\prt#1{$(#1,\le_{#1})$}
\def\eprt#1{(#1,\le_{#1})}
\def\ca#1{$\cal{#1}$}
\def\eca#1{\cal{#1}}
\def\pri{$\le'$}
\def\pris{$<'$}
\def\epri{\le'}
\def\epris{<'}
\def\cas{$\cal{S}$}
\def\ecas{\cal{S}}
\def\cat{$\cal{T}$}
\def\ecat{\cal{T}}
\def\gen#1#2{$\hbox{gen}(#1,#2)$}
\def\egen#1#2{\hbox{gen}(#1,#2)}
\def\gend{$\hbox{gen}(\le_P,S)$}
\def\egend{\hbox{gen}(\le_P,S)}
\def\lad{\{}
\def\rad{\}}
\def\lex{lexicographic }
\def\dom{domain }
\def\wrt{with respect to }
\def\desc{descendable }
\def\nondesc{nondescendable }
\def\arb{arbitrary }

%**A characterization of undirected branching greedoids
% (finite, greedoid minor char undir branch greeds, no antimatroids)
% Wolfgang Schmidt
% http://www.sciencedirect.com.proxy.lib.ohio-state.edu/science/article/pii/0095895688900676
%**Minor characterization of undirected branching greedoids - a short proof
%  (finite, greedoid minor char undir branch greeds, no antimatroids,
% just shorter proof of schmidt's result)
% http://www.sciencedirect.com.proxy.lib.ohio-state.edu/science/article/pii/0012365X9090048M
%**A Characterization of Mixed Branching Greedoids
% (finite, greedoid forb minor char, mixed branching greeds)
% http://www.springerlink.com.proxy.lib.ohio-state.edu/content/kq266k27r2255p74/fulltext.pdf
%**On Paths of Greedoids and a Minor Characterization
% (finite, greedoid minor char interval greeds, mentions antimatroids)
% Erhard Hexel
% http://www.sciencedirect.com.proxy.lib.ohio-state.edu/science/article/pii/S0195669884710250

%korte and lovasz introduce greedoids (where?)
% (finite, compares matroids and antimatroids, basic properties,
%korte lovasz schrader greedoids springer berlin 1991
%diestel infinite matroid thing?
%add forbidden greedoid minor results in lovasz. add references

%http://www.math.osu.edu/~carlson/research/betaadraft.pdf

The following definition of {\em autonomous system} is basic to all that
follows.
An autonomous system is essentially an abstract proof system, though this
is likely not obvious from the definition. We refer the curious reader to
[1], in which this intuition is explained in detail.

\begin{dfn}\label{defausys}
An autonomous system is a set $S$ together with a set $T$ of subsets of
$S$, called autonomous sets, satisfying the following two
conditions:
\begin{itemize}
\item[(i)] $T$ is closed under arbitrary union.
\item[(ii)] (The Order Property) For every autonomous set $A$
there is a total order $\le'$ of $A$ such that every $\le'$
downward closed set is autonomous.
\end{itemize}
\end{dfn}

Note that since the union of autonomous sets is autonomous, there
is a largest autonomous subset $S'$ of $S$. Since all the structure
is contained within $S'$, we assume unless otherwise noted that
$S=S'$. More generally, every set $X$ in an autonomous system has
a largest autonomous subset $A$ under inclusion. We call $A$ the
autonomous part of $X$. We call $X-A$ the nonautonomous part of $X$.

\section{The Canonical Orders}

The canonical orders may be thought of as context dependent orders of
needing. For an autonomous set $A$, we think of $x<_A y$ as saying that
$y$ needs $x$ if one is restricted to only using tools from the set $A$.
If the notion of proof is defined abstractly, then $x<_A y$ means that
$x$ precedes $y$ in every proof whose underlying set is a subset of $A$.
Since the definition of proof is outside the scope of this paper, we
take the following as our definition, though it is normally a proposition.

\begin{dfn} \label{xlessychar}
  Let $A$ be autonomous and let $x,y$ be in $A$.
  Then $x\ere{A}y$ iff every \aut subset
  of $A$ containing $y$ also contains $x$.
\end{dfn}

We need several lemmas.

\begin{lem}\label{minautmaxcan}
Let $S$ be an autonomous system, let $x$ be in $S$, and let $A$ be a minimal
\aut subset of $S$ containing $x$. Then $x$ is a maximum element
of the canonical order $\le_A$.
\end{lem}

\begin{prf}
Suppose not. Then there is $y$ in $A$ such that $y\not<_A x$.
Therefore there is an \aut subset $B$ of $A$ containing $x$
and not $y$. So $A$ is not a minimal \aut set containing $x$,
a contradiction.
\end{prf}

\begin{lem}\label{canordersweakenone}
  Let $A\subseteq B$ be \aut sets and $x,y\in A$. If $x\eres{B}y$
  then $x\eres{A}y$.
\end{lem}

\begin{prf}
  Suppose $x\eres{B}y$. Then
  every \aut subset of $B$ containing $y$ also contains $x$.
  In particular, every \aut subset of $A$ containing $y$ also contains $x$.
  Therefore $ x\eres{A}y$.
\end{prf}

\begin{lem}\label{autdcincan}
  Let $A$ be \aut and $B$ an \aut subset of $A$. Then $B$ is \dc in the
  canonical order \re{A}.
\end{lem}

\begin{prf}
  We have to show $B$ is \re{A}
  downward closed, so let $y\in B$ and $x\eres{A}y$. We must show
  $x\in B$. Since $x\eres{A}y$, we see every \aut subset of $A$
  containing $y$ also contains $x$. In particular, $B$ contains $x$.
\end{prf}

\section{Partial Orders As Autonomous Systems}

\begin{lem}\label{dcsatisfiesorderproperty}
Let $(P,\le)$ be a partial order. Then the set $T$ of downward
closed sets is closed under arbitrary union and satisfies the
order property.
\end{lem}

\begin{prf}
To see $T$ is closed under arbitrary union, let $D_i$ be downward
closed for each $i$ in an index set $I$. We must show $D=\bigcup_{i\in I} D_i$
is downward closed. Let $y\in D$ and $x<y$. Since $y$ is in $D$ then
$y$ is in some $D_i$. Since $D_i$ is downward closed and $x<y$, we see
that $x$ is in $D_i$. Since $D_i\subseteq D$, we see that $x$ is in $D$.

To see that $T$ satisfies the order property, let $D\in T$ be a
downward closed set. We must show there is a total order $\le'$ on $D$
such that every $\le'$ \dc set is in $T$. In other words, every $\le'$
\dc set must be $\le$ downward closed. This means exactly that
$x\le y$ implies $x\le' y$ for all $x,y$ in $P$. Such a total order is
called a linear extension of $\le$.
Linear extensions are well known to exist for every partial order.
The proof is thus complete.
\end{prf}

\begin{cor}\label{posetisausys}
Let $(P,\le)$ be a partial order and let $T$ be the set of $\le$
\dc sets. Then $(P,T)$ is an autonomous system.
\end{cor}

%We now describe the canonical orders in a partial order.

%\begin{lem}\label{posetcanordersdontchange}
%  Let $A$ and $B$ be \aut sets in a partial order $(P,\le)$ and
%  let $x,y\in A\cap B$. Then $x\ere{A}y$ iff $x\ere{B}y$.
%\end{lem}

%\begin{prf}
%Let $T$ be the set of \re{} \dc sets.
%We know that $x\le_A y$ iff every $T$ \aut subset of $A$
%containing $y$ also contains $x$. Since $T$ is the set of
%\re{} \dc sets, this means $x\le_A y$ iff every \re{} \dc
%subset of $A$ containing $y$ also contains $x$. Since $x$
%is in $A$, this is equivalent to saying $x\le y$. Similarly,
%$x\le_B y$ iff $x\le y$. Therefore $x\le_B y$ iff $x\le_A y$
%as claimed.
%\end{prf}

%\begin{cor}\label{posetcanorderofaut}
%  Let $(P,\le)$ be a partial order, $T$ the set of \re{} \dc subsets of $P$,
%  and $A$ in $T$. In the autonomous system $(P,T)$,
%  the canonical order \re{A} of $A$ is the restriction $\ere{P}|A$
%  of the canonical order of $P$.
%\end{cor}

The next lemma is instrumental in proving a useful characterization of
partial orders as proof systems.

\begin{lem}\label{autintersectdc}
  Let $A$ be \aut and $S$ an \arb subset of $A$. Then $S$ is \re{A} \dc
  iff it is the (possibly infinite) intersection of \aut subsets
  of $A$.
\end{lem}

\begin{prf}
Let us first see that if $S_i$ is an autonomous subset of $A$
for all $i\in I$ then $\bigcap_{i\in I}S_i$
  is \re{A} downward closed. Each $S_i$ is by hypothesis an \aut
  subset of $A$ and so is \re{A} \dc by Lemma
  \ref{autdcincan}. Since
  \dc subsets of an \arb \z are closed under intersection, in
  particular so are the \re{A} \dc subsets. It follows that
  $\bigcap_{i\in I}S_i$ is \re{A} downward closed.

  Now, for the nontrivial direction. We must show every \re{A} \dc
  set $S$ can be represented as $\bigcap_{i\in I}S_i$ for some \aut
  subsets $S_i$ of $A$. It is enough
  to show $S$ is the intersection of all \aut subsets of $A$ containing it,
  so let $I$ index all these sets $S_i$.  It is
  obvious $S$ is a subset of the intersection of all \aut subsets of $A$
  containing it, so we have only to show the reverse inclusion.

  So we have to show the intersection $\bigcap_{i\in I}S_i$ is
  contained in $S$, which means we must show every element of
  $\bigcap_{i\in I}S_i$ is also an element of $S$. We show the
  contrapositive, that given $x\in A$, if $x$ is not an element
  of $S$ then $x$ is not an element of $\bigcap_{i\in I}S_i$.

  So take $x$ not in $S$. To show $x$ is not in $\bigcap_{i\in I}S_i$
  is to show there is $i\in I$ \st $x$ is not in $S_i$. Since
  our $S_i$'s are all the \aut subsets of $A$ containing $S$, this
  means we have to give an \aut subset of $A$ containing $S$ but
  not containing $x$. It is enough to give, for each $s\in S$,
  an \aut subset $B$ of $A$ containing $s$ and not $x$. For then
  $\bigcup_{s\in S}B_i$ will be the desired \aut subset of $A$ containing
  $S$ and not $x$, completing the proof.

  So take $s\in S$. How do show there is an \aut subset of $A$ containing
  $s$ and not $x$? If there were no such \aut subset of $A$,
  then every \aut subset of $A$ containing $s$ would also contain $x$,
  and therefore we would have $x\ere{A}s$. Now $s$ is in $S$ and $S$
  is \re{A} \dc by hypothesis, which implies $x$ is in $S$, contrary
  to our choice of $x$ as an element not in $S$. This contradiction
  proves the lemma.
\end{prf}

The following theorem characterizes partial orders in terms of \aut sets.

%undone maybe add to next lemma two other conditions.
% first about finite and descending intersection
% second about descending intersection and can orders never fatten

\begin{thm}\label{autposetchar}
Let $(P,T)$ be an autonomous system given by \aut sets. Then the following are
equivalent:
\begin{itemize}
\item[(i)] $(P,T)$ is a partial order.
\item[(ii)] The $T$ \aut sets are closed under arbitrary intersection.
\end{itemize}
\end{thm}

\begin{prf}
  $\it{(i)}\Rightarrow \it{(ii)}$:\quad
If $(P,T)$ is a partial order then we may take \re{} on $P$ \st the sets
in $T$ are exactly the \re{} \dc sets. Since the \dc sets of a
\z are closed under \arb intersection, we see the \aut sets of
$(P,T)$ are as well.

  $\it{(ii)}\Rightarrow \it{(i)}$:\quad
  For the converse, we assume the $T$ \aut sets are closed under
  \arb intersection. By Lemma
\ref{autintersectdc}, the \re{P} \dc subsets of $P$ are
exactly the intersections of \aut subsets of $P$. Since we are
assuming the \arb intersection of \aut sets is autonomous, this
implies the \re{P} \dc sets are exactly the $T$ \aut sets. Therefore
$(P,T)$ is a \z as claimed.
\end{prf}

%\begin{lem}\label{}
%  Let $(P,T)$ be an \arb autonomous system. Then the following are equivalent.
%  \begin{itemize}
%    \item[(i)]
%      For all \aut sets $A,B$ and for all $x,y$ in $A\cap B$,
%      $x\ere{A}y\Leftrightarrow x\ere{B}y$.
%    \item[(ii)]
%      For all
%      \aut sets $A\subseteq B$ the orders \re{A} and $\ere{B}|A$
%      are the same.
%  \end{itemize}
%\end{lem}

%\begin{prf}
%  $\it{(i)}\Rightarrow \it{(ii)}$ should be obvious,
%so let us prove $\it{(ii)}\Rightarrow \it{(i)}$.
%  We prove the contrapositive instead, so assume $\it{(i)}$ fails. Then
%  there are \aut sets $A,B$ and $x,y\in A\cap B$ \st $x\ere{A}y$
%  and $x\not\ere{B}y$. Now $B\subseteq A\cup B$ and so
%  \re{A\cup B} \weas \re{B}. Since $x\not\ere{B}y$ this implies
%  $x\not\ere{A\cup B}y$. So we have $x\ere{A}y$ and
%  $x\not\ere{A\cup B}y$, and $A$ is a subset of $A\cup B$. Hence
%  the canonical orders \wea from $A$ to $A\cup B$, and $\it{(ii)}$ fails.
%\end{prf}

\section{Deletion, Contraction, Quotients, and Minors}

We now rigorously define the containment relations
for autonomous systems with which our main theorems are stated.
We first define deletions and contractions.

\begin{dfn}
Let $(P,T)$ be an autonomous system and $C$ a subset of $P$. Then the contraction
$P/C$ of $P$ to $P-C$ is defined
as the autonomous system $(P-C,T')$ with domain $P-C$ and set of autonomous sets
$$T'=\{B\subseteq C: B=A-C \hbox{ for some  } A\hbox{ in } T\}$$
\end{dfn}

\begin{dfn}
Let $(P,T)$ be an autonomous system and let $C$ be a subset of $P$.
Then the deletion $P\backslash C$ of $P$ to $P-C$
is defined as the autonomous system
$(C, T')$ with domain $P-C$ and set of autonomous sets
$$T'=\{B\subseteq C: B\in T\}$$
\end{dfn}

Note that we sometimes denote $P/C$ instead by $P|(P-C)$ and refer
to restricting $P$ to $P-C$. Similarly, we sometimes denote $P\backslash C$
by $P.(P-C)$ and refer to dotting to $P-C$. We refer to an autonomous system
obtained from $P$ by a sequence of deletions and contractions (equivalently
a sequence of dottings and restrictions) as a subdot or delecontraction.

We now define homomorphisms and quotients. These notions, together
with deletions and contractions, will allow us to define minors and induced
minors.

\begin{dfn}
An autonomous system homomorphism is a function $f$ from an
autonomous system $P$ to an autonomous system $Q$ such that
$f^{-1}(A)$ is autonomous in $P$ for all autonomous $A\subseteq Q$.
\end{dfn}

\begin{dfn}
Let $(P,T_P)$ and $(Q,T_Q)$ be autonomous systems.
A surjective autonomous system homomorphism $f:(P,Q_P)\to (Q,T_Q)$ is
called a quotient map if $T'\subseteq T_Q$
for all autonomous systems $(Q,T')$ with domain $Q$
such that $f:(P,Q_P)\to (Q,T')$ is a homomorphism.
The autonomous system $(Q,T_Q)$ is then called a quotient of $(P,T_P)$.
\end{dfn}

\begin{dfn}\label{minor}
Let $Q$ and $Q'$ be autonomous systems. We say that $Q'$ is an induced minor of $Q$
if there is a sequence $Q=Q_1,\ldots, Q_n=Q'$ of autonomous systems such that
for each $i$ with $1\le i<n$, one of the following conditions holds:
\begin{itemize}
\item[(i)] $Q_{i+1}=Q_i/C$ for some subset $C$ of $Q_i$.
\item[(ii)] $Q_{i+1}=Q_i\backslash C$
for some  subset $C$ of $Q_i$.
\item[(iii)] $Q_{i+1}$ is a quotient of $Q_i$.
\end{itemize}
\end{dfn}

If we replace the third condition with the condition that $Q_{i+1}$
is simply a homomorphic image of $Q_i$, we get the notion of autonomous system
minor. The names minor and induced minor are chosen for good reason.
Though technical and outside the scope of this paper, roughly speaking,
it can be shown that considering each graph as a family of autonomous
systems, graph  minor and induced minor correspond to autonomous system
minor and induced minor, respectively.
Allowing homomorphic images that are not quotient maps corresponds to
graph edge deletion.

\section{Autonomous System Join}

Joins of autonomous systems allow us to prove that when an equivalence
relation is homomorphism induced, it is in fact induced by a quotient
map.
We recall that given a partial order $(Z,\le)$ and points $x,y$ in $Z$,
the join $x\lor y$ of $x$ and $y$ is the least upper bound of $x$ and
$y$ if one exists. Otherwise $x\lor y$ is undefined. More generally,
if $S$ is a nonempty subset of $Z$, then $\bigvee_{x\in S}x$ is a
least upper bound of $S$ if one exists and is otherwise undefined.

Given a set $\{(P_i,T_i)\}_{i\in I}$  of autonomous systems, we let
$P=\bigcup_{i\in I}P_i$ and consider the autonomous systems
$(P,T_i)$ for $i$ in $I$. We let $(P,T_i)\le (P,T_j)$ if $T_i\subseteq T_j$.
With this definition of $\le$, we may then speak of the least upper bound,
or join, of a nonempty set of autonomous systems.
The next lemma shows that the join of every nonempty
set of autonomous systems exists.

\begin{lem}
Given a set of autonomous systems $(P_i,T_i)$ for $i$ in a nonempty index
set $I$, the join $\bigvee_{i\in I}(P_i,T_i)$ exists. Specifically, it is
the autonomous system $(P,T)$ on $P=\bigcup_{i\in I}P_i$, where
$T$ is the closure of $\bigcup_{i\in I}T_i$ under arbitrary union.
\end{lem}

\begin{prf}
First, we consider each autonomous system $(P_i,T_i)$ instead as the autonomous system $(P,T_i)$.
As in the statement of the lemma,
we let $T$ be the set of subsets of $P$ of the form $\bigcup_{i\in I} A_i$,
where each $A_i$ is a (possibly empty) set in $T_i$.  From
the definition of $T$ and the fact that autonomous sets in an autonomous
system are closed under arbitrary union, it is immediate that if $(P,T)$
is an autonomous system, then it is in fact the least upper bound as
required. To show that $(P,T)$ is an autonomous system, we must show $T$ is closed
under arbitrary union and satisfies the order property.

Closure under arbitrary union is immediate by definition of $T$. To show
the order property, we must show every set $A$ in $T$ can be totally ordered
such that each downward closed set is also in $T$. We know that
$A=\bigcup_{i\in I}A_i$ for sets $A_i$ in $T_i$. Since $(P,T_i)$ is an
autonomous system for each $i$ in $I$, we see that each $T_i$ satisfies the order
property. We may therefore totally order each $A_i$ as $\le_i$ such that
each $\le_i$ downward closed set is autonomous.

Choose a well order $\le_w$ on the set $I$.
We  define a total order $\le$ on $A$ as follows.
For each $x$ in $A$, let $r_x$ be the $\le_w$ least element of $I$ such
that $x$ is in $A_i$. If $r_x<_w r_y$ then let $x<y$.
If $r_x=r_y=i$ then let $x\le y$ iff $x\le_i y$.
The reader may check that this is a well defined total order on $A$ such
that each $\le$ downward closed set $D$ is the union of $\le_i$ downward
closed sets $D_i$ for each $i$. Since each $D_i$ is $\le_i$ downward
closed, it follows that $D_i$ is $T_i$ autonomous. Since $T_i\subseteq T$,
we see that $D_i$ is $T$ autonomous. Since $T$ is closed under arbitrary
union, it follows that
$D=\bigcup_{i\in I}D_i$
is $T$ autonomous. This completes the proof.
\end{prf}

\begin{dfn}
Let $f:P\to Q$ be an autonomous system homomorphism. We say the equivalence relation
$\sim$ on $P$ such that $x\sim y$ iff $f(x)=f(y)$ is induced by $f$. We call
an equivalence relation on $P$ homomorphism induced if it is induced by some
surjective
homomorphism $f:P\to Q$ such that $Q$ has at least one nonempty autonomous
set.
\end{dfn}

The requirement that $Q$ has at least one nonempty autonomous set is given
because otherwise every equivalence relation would vacuously be homomorphism
induced by a map to an autonomous system with only the empty set autonomous.

\begin{lem}
If $f:P\to (Q,T_i)$ is an autonomous system homomorphism for each $i$ in a
nonempty index
set $I$, then $f:P\to \bigvee_{i\in I} (Q,T_i)$ is a homomorphism.
\end{lem}

\begin{prf}
We must show the inverse image of every $\bigvee_{i\in I} (Q,T_i)$
autonomous set is $P$ autonomous, so choose such a set $A$.
Then $A=\bigcup_{i\in I} A_i$ for some sets $A_i$ in $T_i$. For each $i$,
$f:P\to (Q,T_i)$ is a homomorphism and therefore the set
$f^{-1}(A_i)$ is $P$ autonomous. Since the $P$ autonomous sets are closed
under arbitrary union, we see that
$$f^{-1}(A)=f^{-1}(\bigcup_{i\in I}A_i)=\bigcup_{i\in I}f^{-1}(A_i)$$
is $P$ autonomous as needed.
\end{prf}

\begin{lem}
If $P$ is an autonomous system, $Q$ is a set, $f$ a function from $P$ to $Q$,
and $\{(Q,T_i):i\in I\}$ is the set of all autonomous systems on $Q$ such that
$f:P\to (Q,T_i)$ is an autonomous system homomorphism,
then $f:P\to \bigvee_{i\in I} (Q,T_i)$ is a quotient map.
\end{lem}

\begin{prf}
Let $(Q,T)=\bigvee_{i\in I} (Q,T_i)$.
It is only to show that if $f:P\to (Q,T')$ is an autonomous system homomorphism
then $T'\subset T$. This is immediate from the definition of join
and $(Q,T)$.
\end{prf}

\begin{cor}
If an equivalence relation on an autonomous system is homomorphism induced,
then it is induced by a quotient map.
\end{cor}

If $f:P\to (Q,T_1)$ and $f:P\to (Q,T_2)$ are quotient maps, then
$T_1\subseteq T_2\subseteq T_1$ by definition of quotient map,
so $T_1=T_2$. Thus the quotient map of the previous corollary is
unique. We may thus refer to the quotient $P/\sim$ for any homomorphism
induced equivalance relation on $P$. The reader should note that
$P/\sim$ for an equivalance relation $\sim$ and $P/X$ for a subset $X$
of $P$ are distinct notions.

\section{Strong Aut Descendability}

Our main theorems will be stated for the class of strong aut descendable
autonomous systems, which includes both finite autonomous systems and
arbitrary partial orders. The reader who is content to consider finite
systems may skip this section and insert ``finite'' everywhere he or
she reads strong aut descendable.

\begin{dfn}\label{strongaut}
An autonomous system is strong aut descendable if the intersection of every
chain of \aut sets under inclusion is also autonomous.
\end{dfn}

The condition is meant to allow infinite autonomous systems while excluding
certain pathologies that prevent many statements from being true for arbitrary
autonomous systems. The simplest example is an autonomous system on a set
$S\cup x$ such that $S$ is an infinite set not containing $x$, every subset
of $S$ is autonomous, and every cofinite subset of $S\cup x$ is autonomous.
The reader may show that this is an autonomous system. It is somewhat of an
all purpose counterexample in the sense that for each statement we make only
assuming strong aut descendability, the counterexample to the general
statement is very similar in spirit to that just described.

The proof of the following lemma may be taken as a simple exercise in Zorn's
Lemma.

\begin{lem}
If $A$ is an autonomous set in a strong aut descendable autonomous system
and $x$ is a point in $A$, then there is a minimal autonomous subset $B$
of $A$ containing $x$.
\end{lem}

We need to show that strong aut descendability is preserved under taking
subdots.

%\begin{lem}
%An autonomous system $P$ is strong aut descendable iff
%$\bigcap_{i<\alpha} A_i$ is autonomous for each ordinal
%$\alpha$ and for each $\alpha$ sequence $\{A_i\}_{i<\alpha}$ of autonomous
%subsets of $P$ such that $A_j\subset A_i$ for each $i<j<\alpha$.
%\end{lem}

%\begin{prf}
%Of course if $P$ is strong aut descendable, then since the intersection
%of every chain of autonomous sets under inclusion is autonomous, in
%particular this condition holds for descending chains indexed by an
%ordinal.

%For the converse, suppose
%$\bigcap_{i<\alpha} A_i$ is autonomous for each ordinal
%$\alpha$ and for each $\alpha$ sequence $\{A_i\}_{i<\alpha}$ of autonomous
%subsets of $P$ such that $A_j\subset A_i$ for each $i<j<\alpha$.
%Let $C$ be an arbitrary chain of autonomous sets under inclusion. The
%sets in $C$ may be enumerated as $\{B_i\}_{i<\beta}$ for some ordinal
%$\beta$. Now let
%$$A_i=\bigcap_{j\le i}B_j.$$
%The sequence of $A_i$'s is decreasing under inclusion. Taking a subsequence
%if necessary, we may assume without loss of generality that
%$\{A_i\}_{i<\alpha}$ is strictly decreasing under inclusion. Since the
%intersection of such families is autonomous by hypothesis, and since
%$$\bigcap_{i<\alpha}A_i=\bigcap_{i<\beta}B_i,$$
%which is the intersection of the chain $C$ of sets by construction,
%the proof is complete.
%\end{prf}

\begin{lem}\label{subdot_strong_aut_desc}
Strong aut descendability is preserved under taking subdots.
\end{lem}

\begin{prf}
We must show that strong aut descendability is preserved under dotting
and restricting. Let $P$ be an autonomous system. Dotting yields an
autonomous system of the form $P.A$, where $A\subseteq P$ is autonomous.
Since $P$ is strong aut descendable, the intersection of every chain of
subsets of $P$
is autonomous in $P$. In particular, the intersection of every chain
of subsets of $A$
is a subset of $A$ that is autonomous in $P$. The intersection is therefore
autonomous in $P.A$, proving that $P.A$ is strong aut descendable.

To show that strong aut descendability is preserved under restriction,
let $X$ be an arbitrary subset of $P$ and consider a chain $(I,\le)$ of
$P|X$ autonomous sets $\{X_i\}_{i\in I}$ under inclusion such that
$i<j$ iff $X_i\subset X_j$. By definition of restriction, each
$X_i$ has the form $A_i\cap X$ for some some $P$ autonomous set $A_i$.
Given $i$, let $$B_i=\bigcup_{j\le i}A_j.$$
Then each $B_i$ is $P$ autonomous as well since the autonomous sets are
closed under arbitrary union. Moreover, $\{B_i\}_{i\in I}$ is a chain
under inclusion so that $\bigcap_{i\in I}B_i$ is autonomous by strong
aut descendability of $P$. Therefore

$$X\cap\left(\bigcap_{i\in I}B_i\right)
=X\cap\left(\bigcap_{i\in I}\bigcup_{j\le i}A_j\right)
=X\cap\left(\bigcap_{i\in I}A_i\right)
=\bigcap_{i\in I}(X\cap A_i)
=\bigcap_{i\in I}X_i$$
is $P|X$ autonomous as needed.
\end{prf}

\section{The Main Theorems}

The author showed in his doctoral thesis that partial orders comprise an
induced minor closed class of autonomous systems. We omit the somewhat
technical proof.

\begin{lem}
If $(P,\le)$ is a partial order and $(Q,T_Q)$ is an induced minor of $P$,
then $Q$ is a partial order.
\end{lem}

Knowing that partial orders form an induced minor closed class, it is natural
to seek a forbidden induced minor characterization of this class. The main step
in proving such a theorem is Lemma \ref{tcsubdot}. It is stated in terms of the
autonomous system $P_3$. In general, for $n\ge 1$, the autonomous system
$P_n$ has the vertex set $\{v_1,\ldots,v_n\}$ of an $n$ point path as its
underlying set, with a subset of $P_n$ autonomous iff it has the form
$\{v_1,\ldots,v_i\}\cup\{v_k,\ldots,v_n\}$ with $0\le i\le n$ and
$1\le k \le n+1$. In other words, the autonomous sets are the unions of
paths starting at the endpoints and empty paths.

\begin{lem} \label{tcsubdot}
Let $P$ be a strong aut descendable autonomous system with autonomous sets $A$
and $B$ such that $A\cap B$ is not autonomous. Then $P$ contains
$P_3$ as a subdot.
\end{lem}

\begin{prf}
By dotting to $A\cup B$ if necessary, we may assume $P=A\cup B$.
Since $A\cap B$ is not autonomous, we may choose $x$ in the nonautonomous
part of $A\cap B$. Let $A'=(A-B)\cup\{x\}$, $B'=(B-A)\cup\{x\}$,
let $S=(A'\cup B'\cup\{x\})$, and let $P'=P|S$. Since $A'=A\cap S$
and $B'=B\cap S$, we see by definition of $P'$ that $A'$ and $B'$
are autonomous in $P'$. If $\{x\}=A'\cap B'$ is autonomous in $P'$
then by definition of $P'$ as a restriction, $P$ must contain an
autonomous subset of $A\cap B$ containing $x$, contrary to choice
of $x$ in the nonautonomous part. This contradiction shows
$\{x\}=A'\cap B'$ is not autonomous in $P'$.

Since $P'$ is strong aut descendable by Lemma \ref{subdot_strong_aut_desc},
we may choose a minimal $P'$
autonomous subset $A''$ of $A'$ containing $x$. Similarly, we may
choose a minimal $P'$ autonomous subset $B''$ of $B'$ containing $x$.
Note that $A''\cup B''$ is $P'$ autonomous. Let $P''=P.(A''\cup B'')$. Then
$A''$ and $B''$ are $P''$ autonomous sets such that $A''\cap B''=\{x\}$
is not $P''$ autonomous.

Note that since $A''$ is a minimal $P'$ autonomous subset of $A'$
containing $x$, it follows that $x$ is a $\le_{A''}$ maximum element
of $A''$ in the autonomous system $P'$. Since $P'$ and $P''$ have the same
autonomous subsets of $A''$, it follows that $x$ is a $\le_{A''}$
maximum element of $A''$ in the autonomous system $P''$ as well. Since $\{x\}$
is not autonomous in $P''$, it follows that there is $a$ in $A''$
such that $a<_{A''}x$ in the autonomous system $P''$. Similarly,
there is $b$ in $B''$ such that $b<_{B''}x$ in the autonomous system $P''$.

Let $P'''=P''|\{a,b,x\}$. Simple checking of the autonomous sets
of $P'''$ shows that $P'''$ is isomorphic to $P_3$. Obviously,
$P'''$ is a subdot of $P$.
\end{prf}

\begin{lem}\label{sad_nonposet_two_nonintersect}
If a strong aut descendable autonomous system $P$ is not a partial order,
then there are two autonomous sets whose intersection is not autonomous.
\end{lem}

\begin{prf}
Since $P$ is not a partial order,
it follows from Theorem \ref{autposetchar}  that
the autonomous sets are not closed under arbitrary intersection.
Let $\{A_i\}_{i<\lambda}$ be a family of autonomous sets in $P$
whose intersection is not autonomous, for some finite or infinite
cardinal $\lambda$. For each $i<\lambda$, let
$$B_i=\bigcap_{j\le i}A_j.$$
Note that the $B_i$'s comprise a descending chain of sets under inclusion whose
intersection is not autonomous. Therefore some $B_i$ is not autonomous.
Choose the least such $i$. Then $i$ can not be a limit ordinal, for then
$$B_i=\bigcap_{j\le i}A_j=\bigcap_{j< i}A_j,$$
so $B_i$ would be the intersection of a chain of autonomous sets and
therefore autonomous.

So the least such $i$ must be a successor ordinal. That is, $i=j+1$.
Therefore $B_i=B_j\cap A_{j+1}$. Since $B_j$ and $A_{j+1}$ are autonomous
but $B_i$ is not, the proof is complete.
\end{prf}

%proofread once til here
\begin{cor} \label{posetiffnotwochoicesubdot}
A strong aut \desc autonomous system is a partial order iff it has no $P_3$ subdot.
\end{cor}

\begin{prf}
We know the subdot of a partial order is a partial order, which is therefore
not $P_3$ as $P_3$ contains the autonomous sets $\{v_1,v_2\}$ and $\{v_3,v_2\}$
whose intersection is not autonomous.
Inversely, if a strong aut \desc autonomous system $P$
is not a partial order, then by Lemma \ref{sad_nonposet_two_nonintersect}
there are \aut sets $A$ and $B$ \st $A\cap B$
is not autonomous.
By Lemma \ref{tcsubdot}, it follows the autonomous system has a $P_3$ subdot.
\end{prf}

\begin{cor}
A strong aut \desc autonomous system is a partial order iff it has no $P_3$
induced minor.
\end{cor}

\begin{prf}
Immediate from the fact that a subdot is a minor and the fact
that the pure contraction of a partial order is a partial order.
\end{prf}

Consider the autonomous system $P_4$ with vertices $a,x,y,b$ in that order
in the path. Note that
$x<_{\{a,x,y\}}y$ and $y<_{\{b,y,x\}}x$. In fact, $P_4$ is the simplest
autonomous system exhibiting such behavior in the sense that every strong aut \desc
autonomous system containing points $x$ and $y$ and sets $A$ and $B$ such that
$x<_A y$ and $y<_B x$ contains a $P_4$ induced minor. We prove this claim now.

The proof that follows is very similar in spirit to the proof
of Lemma \ref{tcsubdot}, which
suggests that the following theorem and
Corollary \ref{posetiffnotwochoicesubdot} may be combined.
In fact Corollary \ref{posetiffnotwochoicesubdot} and the following
theorem may be seen as corollaries of the same result, but that result
is far more technical to state and prove, so we treat each separately.

\begin{thm}\label{pfourminor}
Let $P$ be a strong aut \desc autonomous system with \aut sets $A$ and $B$ and
points $x$ and $y$ such that $x<_A y$ and $y<_B x$. Then $P$ contains
a $P_4$ induced minor.
\end{thm}

\begin{prf}
By dotting to $A\cup B$ if necessary, we may assume $P=A\cup B$.
Given such $A$, $B$, $x$, and $y$,
let $A'=(A-B)\cup\{x,y\}$, $B'=(B-A)\cup\{x,y\}$,
let $S=(A'\cup B'\cup\{x,y\})$, and let $P'=P|S$. Since $A'=A\cap S$
and $B'=B\cap S$, we see by definition of $P'$ that $A'$ and $B'$
are both autonomous sets in $P'$ containing $x$ and $y$. If there is a
$P'$ \aut subset $C$ of $A'$ containing $y$ and not $x$ then
$C=D\cap S$ for some $P$ autonomous set $D$. But then $D\subseteq A$,
contrary to the fact that there is no $P$ \aut subset of $A$ containing
$y$ and not $x$. This contradiction shows that every $P'$ \aut subset
of $A'$ containing $y$ also contains $x$. Therefore $x<_{A'}y$ in
the autonomous system $P'$. Similarly, $y<_{B'}x$ in the autonomous system $P'$.

Since $P'$ is strong aut descendable by Lemma \ref{subdot_strong_aut_desc},
we may choose a minimal $P'$
autonomous subset $A''$ of $A'$ containing $y$. Similarly, we may
choose a minimal $P'$ autonomous subset $B''$ of $B'$ containing $x$.
$A''\cup B''$ is $P'$ autonomous. Let $P''=P.(A''\cup B'')$. Then
$A''$ and $B''$ are $P''$ autonomous sets such that $A''\cap B''=\{x,y\}$.

Note that since $A''$ is a minimal $P'$ autonomous subset of $A'$
containing $y$, it follows that $y$ is a $\le_{A''}$ maximum element
of $A''$ in the autonomous system $P'$. In particular, $x<_{A''}y$ in the autonomous system $P'$.
Since $P'$ and $P''$ contain the same \aut subsets of $A''$, we see that
$x<_{A''}y$ in the autonomous system $P''$ as well.
Similarly, $y<_{B''}x$ in the autonomous system $P''$.
Since $x<_{A''}y$ in $P''$, we see that $y$ is
not a $P''$ axiom. Therefore $B''-\{x,y\}$ is nonempty. Similarly,
$A''-\{x,y\}$ is nonempty.

The reader may check that the partition of $P''$ with cells
$A''-\{x,y\}$, $B''-\{x,y\}$, $\{x\}$, and $\{y\}$ is homomorphism
induced. Let $P'''$ be the pure contraction of $P''$ with respect
to this partition. The reader may check that this autonomous system is isomorphic
to $P_4$.
\end{prf}

We note that while our forbidden induced minor theorem for partial orders
can also be seen as a forbidden subdot theorem, the same is not
true for the theorem just stated. Let $P$ be the six point autonomous system
on $\{a_1,a_2,x,y,b_1,b_2\}$
whose autonomous sets are the (possibly empty) unions of the sets
$\{a_1\}$, $\{a_2\}$, $\{b_1\}$, $\{b_2\}$, $\{a_1,x\}$, $
\{a_2,x\}$, $\{b_1,y\}$, $\{b_2,y\}$, $A:=\{a_1,a_2,x,y\}$,
and $B:=\{b_1,b_2,y,x\}$.
The reader may show that this is an autonomous system for which
$x<_A y$ and $y<_B x$, yet $P$ has no $P_4$ subdot.

%\section{Future Directions}

%weaken the strong aut desc hyp
%Note there's one guaranteed induced minor. There is a finite set of guaranteed
%induced minors in a much more general situation, and they can in fact be
%described,
%but doing so is significantly more involved and will be the subject of a
%future paper
%is there a finite set of guaranteed induced minors for unidirectionality

%\bibliographystyle{plain}
%bibliography{mybib}

\,

\,

\,

\,

\,

[1] Christian Altomare. {\it Degree Sequences, Forcibly Chordal Graphs and Combinatorial Proof Systems.}
PhD thesis. The Ohio State University, December 2009. http://etd.ohiolink.edu

[2] Timothy J. Carlson. {\it Ranked $\beta$-Calculi I: Framework.} (Draft version.)
http://www.math.osu.edu/~carlson.6/research/betaadraft.pdf

[3] Brenda L Dietrich. {\it Matroids and antimatroids - a survey.} Discrete Mathematics {\bf 78} (3) 
(1989) pp. 223--237

[4] Robert P. Dilworth. {\it Lattices with unique irreducible decompositions.} Annals of Mathematics {\bf 41} (4)
(1940) pp. 771--777

[5] Paul H. Edelman and Robert E. Jamison. {\it The theory of convex geometries.} Geometriae Dedicata {\bf 19} (3)
(1985) pp. 247--270

[6] Kenji Kashiwabara and Masataka Nakamura. {\it A characterization of the edge-shelling convex geometries of trees.}
http://arxiv.org/abs/0908.3456
% (finite, two antimatroid minor chars, not same as yours,
%  edge shelling convex geometry, ones with all stems size two)

[7] Kenji Kashiwabara and Masataka Nakamura. {\it Characterizations of the convex geometries arising from the double shellings of posets.} Discrete Mathematics {\bf 310} (15-16) (2010) pp. 2100-2112
% (finite, all antimatroid, forb minor char double poset shelling)
% http://www.sciencedirect.com.proxy.lib.ohio-state.edu/science/article/pii/S0012365X1000110X

[8] Masataka Nakamura. {\it Excluded-minor characterizations of antimatroids arisen from posets and graph searches.}
Discrete Applied Mathematics {\bf 129} (2-3) (2003) pp. 487--498
% (finite, all antimatroid, two forb minor chars, shelling and
% node search graph digraph, maybe three chars recheck)

[9]  Masataka Nakamura and Yoshio Okamoto. 
{\it The forbidden minor characterization of line-search antimatroids of rooted digraphs.}
Discrete Applied Mathematics {\bf 131} (2) (2003) pp. 523-533
%  (finite, all antimatroid, forb minor rooted digraph line search)

[10] Nathalie Wahl. {\it Antimatroids of finite character.} Journal of Geometry {\bf 70} (1-2) (2001) pp. 168-175

\end{document}